\documentclass[12pt,dvipsnames,a4paper]{article}
\usepackage[affil-it]{authblk}
\usepackage[utf8]{inputenc}
\usepackage{amsmath, amssymb, mathrsfs, amsthm, amsfonts}
\usepackage{latexsym}
\usepackage{hyperref}
\usepackage{a4wide}
\usepackage{color,xcolor}
\usepackage{graphicx}
\usepackage{fancyhdr}

\usepackage{natbib}

\newtheorem{theorem}{Theorem}[section]

\newtheorem{proposition}[theorem]{Proposition}

\theoremstyle{definition}
\newtheorem{definition}[theorem]{Definition}
\theoremstyle{remark}
\newtheorem*{example}{Example}

\author{Hery Randriamaro
	\thanks{This research was funded by my mother \\
		Lot II B 32 bis Faravohitra, 101 Antananarivo, Madagascar \\
		e-mail: \texttt{hery.randriamaro@outlook.com}}}

\title{Face Counting for Topological Hyperplane Arrangements}

\author{Hery Randriamaro}

\begin{document}

\thispagestyle{empty}

\begin{center}
{\color{MidnightBlue} \bfseries \Huge Face Counting for Topological Hyperplane Arrangements}
\end{center}

\vspace*{20pt}

\noindent \textbf{\Large Hery Randriamaro} 

\vspace*{10pt}

\noindent \textsc{Institut für Mathematik, Universität Kassel,\\ Heinrich-Plett-Straße 40,\\ 34132 Kassel, Germany} \\ \texttt{hery.randriamaro@mathematik.uni-kassel.de} \\ The author was supported by the Alexander von Humboldt Foundation

\vspace*{20pt}

\noindent \textsc{\large Abstract} 

\noindent \noindent Determining the number of pieces after cutting a cake is a classical problem. \cite{roberts1887figures} provided an exact solution by computing the number of chambers contained in a plane cut by lines. About 88 years later, \cite{zaslavsky1975facing} even computed the $f$-polynomial of a hyperplane arrangement, and consequently deduced the number of chambers of that latter. Recently, \cite{forge2009division} introduced the more general structure of topological hyperplane arrangements. This article computes the $f$-polynomial of such arrangements when they are transsective, and therefore deduces their number of chambers.

\vspace*{10pt}


\vspace*{10pt}


\vspace*{10pt}

\noindent {\color{MidnightBlue} \rule{\linewidth}{2pt}}

\vspace*{10pt}

\section{Introduction}

\noindent A classical basic problem was to determine the number of pieces obtained by cutting a cake $d$ times. Deeper study of that problem has probably its origin in the article of \cite{steiner1826einige} who computed the maximal number of chambers contained in a plane cut by several sets of parallel lines pointing in different directions. \cite{roberts1887figures} fixed that problem by showing that $\displaystyle 1 + d + \binom{d}{2} + \sum_{i=1}^k n_k \binom{k-1}{2} - \sum_{j=1}^p \binom{l_j}{2}$ is the number of chambers contained in a plane cut by $d$ lines, where $n_k$ is the number of $k$-fold intersection points for $k \geq 3$, and $p$ is the number of families of parallel lines containing respectively $l_1, \dots, l_p$ lines with $l_j \geq 2$. As mentioned in the book of \cite{dimca2017hyperplane} for instance, Schläfli extended that problem to the Euclidean space $\mathbb{R}^n$, and published in 1901 that the number of chambers in $\mathbb{R}^n$ partitioned by $d$ hyperplanes is smaller that $\displaystyle \sum_{i=0}^n \binom{d}{i}$. That extended problem was, that time, solved by \cite{zaslavsky1975facing}. He precisely expressed the $f$-polynomial of a hyperplane arrangement $\mathcal{A}$ by means of its Möbius polynomial, and deduced that its number of chambers is $\displaystyle \sum_{X \in L(\mathcal{A})} (-1)^{\mathrm{rank}\, X} \mu(\mathbb{R}^n,X)$, where $L(\mathcal{A})$ is the flat set of $\mathcal{A}$ and $\mu$ the Möbius function. In an independent work, \cite{alexanderson1981arrangements} obtained the $f$-polynomial of a plane arrangement in a space. More recently, \cite{pakula2003pseudosphere} computed the number of chambers of pseudosphere arrangements. Note that pseudosphere arrangements are topologically equivalent to pseudohyperplane arrangements as one can read in the article of \cite{deshpande2016arrangements} for example.
 
\smallskip 

\noindent This article considers the more general case of topological hyperplane arrangements, or topoplane arrangements, introduced by \cite{forge2009division}. Transsective topoplane arrangements are even generalizations for pseudohyperplane arrangements that are known to be topological models for oriented matroids, like stated in the book of \cite{bjorner1999oriented}. \emph{This article determines the $f$-polynomial of a transsective topoplane arrangement $\mathscr{A}$ in a topological ball $T$, and deduces that $\displaystyle \sum_{X \in L(\mathscr{A})} (-1)^{\mathrm{rank}\, X} \mu(T,X)$ is its number of chambers, where $L(\mathscr{A})$ is the flat set of $\mathscr{A}$.} 
 
\smallskip
 
\noindent In neighboring contexts, \cite{dumitrescu2019new} established that the number of nonisomorphic simple arrangements of $n$ pseudolines is bigger that $2^{cn^2 - O(n \ln n)}$ for some constant $c > 0.2083$, while \cite{felsner2020arrangements} studied the circularizability of pseudocircle arrangements. Recall that in the Euclidean space $\mathbb{R}^n$, an $n$-ball of radius $r$ and center $x$  is the set of all points of distance less than $r$ from $x$, a topological $n$-ball is any subset which is homeomorphic to an $n$-ball, and an $n$-manifold is a subset with the property that each point has a neighborhood that is homeomorphic to an $n$-ball. Topological $n$-balls are important as building blocks of CW-complexes. However, they are not flexible enough to investigate topological properties of topoplane arrangements. More abstract objects, named deformed $n$-balls, must consequently be introduced in Section~\ref{DeBa}. The study topoplane arrangements really begins in Section~\ref{DeAr}. We namely fix the conjecture of \cite{forge2009division}, mentioned in the introduction of their article, stating that solidity can be proved from the definition of a topoplane arrangement. Then, we prove that every chamber of a transsective topoplane arrangement is a deformed ball. These results allow us to compute the $f$-polynomial of a transsective topoplane arrangement in Section~\ref{fPo}, and to deduce its number of chambers.

\section{Deformed Balls} \label{DeBa}

\noindent This article uses the notations $[k] := \{1, 2, \dots, k\}$ for a positive integer $k$, and $\mathbb{N}_0$ for the set of nonnegative integers. Deformed balls, deformed ball complexes, as well as the Euler characteristic of a deformed ball complex are defined in this section.

\begin{definition}
Let $n$ be a nonnegative integer. A \emph{deformed $n$-ball} is an $n$-manifold $X$ in $\mathbb{R}^n$ that is path connected, and such that the homotopy group $\pi_k(X, x_0)$ is trivial for each positive integer $k$ and a distinguished point $x_0$ of $X$.  
\end{definition}

\begin{definition}
Let $X$ be a deformed $n$-ball, and $Y$ a deformed $m$-ball such that $n>m$ and $X \cap Y = \varnothing$. The sets $X$ and $Y$ can be glued together if the boundary $\partial X$ of $X$ contains $Y$. The set obtained from gluing $Y$ onto $X$ is the path connected space $X \sqcup Y$. 
\end{definition}

\paragraph{Recursive Construction of a System of Deformed Balls.} We begin with a system $\big(X_1,\,\{X_1\}\big)$, where $X_1$ is a deformed $n$-ball.
\begin{itemize}
\item Let $X_2$ be a deformed $m$-ball such that $X_2$ can be glued onto $X_1$, if $n > m$, or $X_1$ can be glued onto $X_2$, if $n < m$. We get the extended system $\big(X_1 \sqcup X_2,\, \{X_1, X_2\}\big)$.
\item Suppose now that we have a positive integer $k$, and a system $\big(X,\,\{X_i\}_{i \in [k]}\big)$, where $\displaystyle X = \bigsqcup_{i \in [k]} X_i$ was obtained by gluing together the deformed balls $X_1, \dots, X_k$. This system can be extended with another deformed ball $X_{k+1}$ if
\begin{itemize}
\item $X \cap X_{k+1} = \varnothing$,
\item there exists $i \in [k]$ such that $X_i$ and $X_{k+1}$ can be glued together,
\item if $I$ is the subset of $[k]$ such that $X_i$ and $X_{k+1}$ can be glued together for each $i \in I$, then $\displaystyle \bigsqcup_{i \in I} X_i$ is path connected.
\end{itemize} 
We then obtain a new system $\big(X \sqcup X_{k+1},\, \{X_i\}_{i \in [k+1]}\big)$ of deformed balls.
\end{itemize}

 \begin{definition}
 A topological space $X$ is a \emph{deformed ball complex} is there exist a positive integer $k$, and a set $\{X_i\}_{i \in [k]}$ of deformed balls such that $\displaystyle X = \bigsqcup_{i \in [k]} X_i$ and $\big(X,\,\{X_i\}_{i \in [k]}\big)$ is a system of deformed balls.
 \end{definition}

\noindent For a CW complex $X$, the Euler characteristic $\chi(X)$ is the alternating sum $\displaystyle \sum_{n \in \mathbb{N}_0} (-1)^n c_n$, where $c_n$ is the number of topological $n$-balls of $X$. We need to generalize that definition deformed ball complexes.

\begin{definition}
Let $k$ be a positive integer, and $\big(X,\,\{X_i\}_{i \in [k]}\big)$ a system of deformed balls. The \emph{Euler characteristic} of the deformed ball complex $X$ is $$\chi(X) := \sum_{n \in \mathbb{N}_0} (-1)^n c_n,$$ where $c_n$ is the number of deformed $n$-balls in $\{X_i\}_{i \in [k]}$.
\end{definition}

\begin{example}
In the left part of Figure~\ref{Comp}, we have a deformed ball complex composed by the deformed $0$-ball, $1$-ball, and $3$-ball represented in the right part of Figure~\ref{Comp}. Its Euler characteristic is $(-1)^0 + (-1)^1 + (-1)^2 = 1$.
	
\bigskip
	
\begin{figure}[h]
\centering
\includegraphics[scale=0.15]{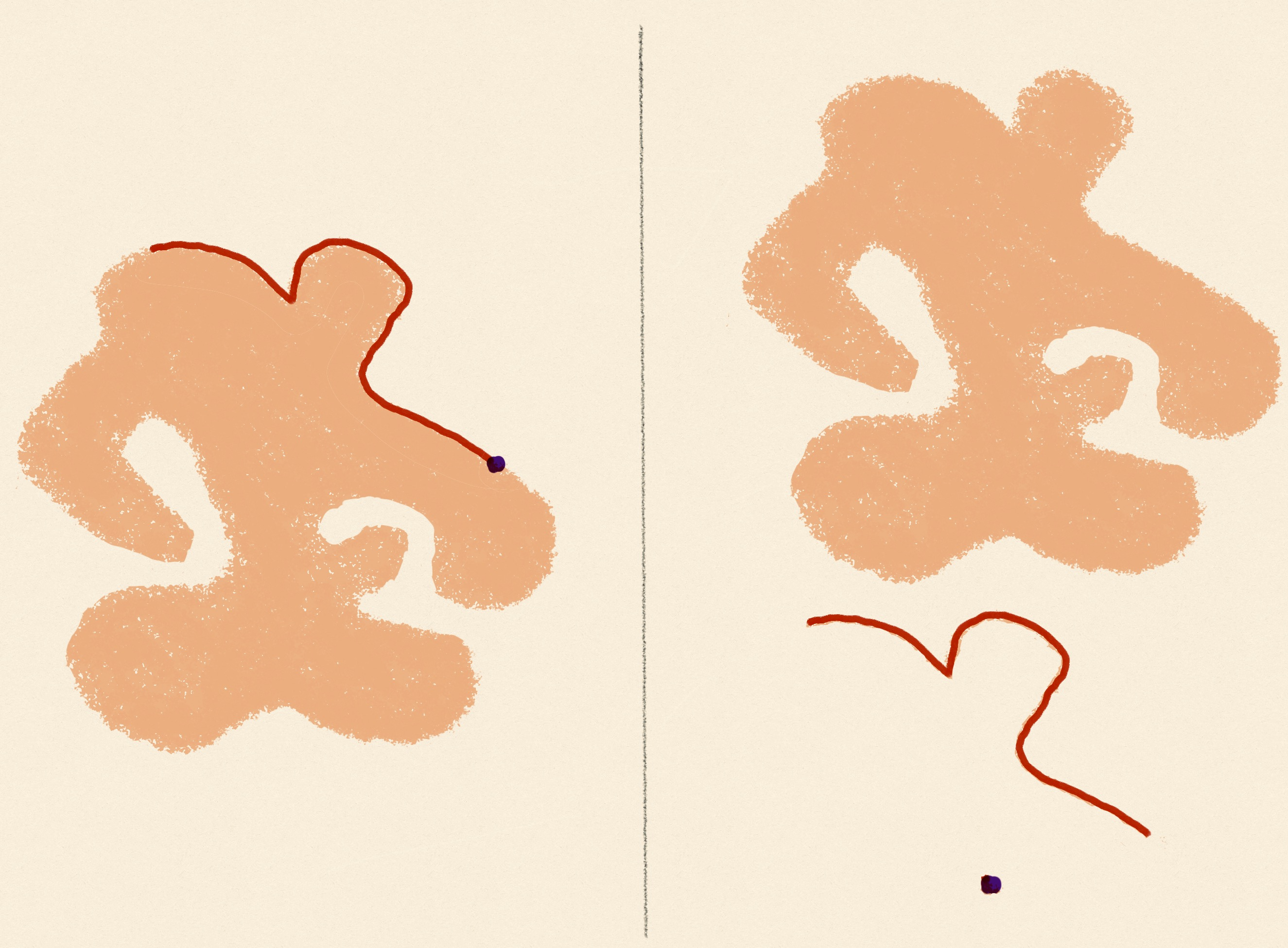}
\caption{A Complex formed by three Deformed Balls}
\label{Comp}
\end{figure}
\end{example}

\section{Topoplane Arrangements}  \label{DeAr}

\noindent This section is devoted to topoplane arrangements introduced by \cite{forge2009division}. Transsective topoplane arrangements are particularly for interest for us. We fix in Proposition~\ref{PrTra} the conjecture mentioned in the introduction of the article of \cite{forge2009division}, stating that every restriction of a transsective topoplane arrangement is a transsective topoplane arrangement. Afterwards, we prove in Proposition~\ref{PrDef} that every face of a transsective topoplane arrangement is a deformed ball.

\begin{definition}
Let $n$ be a positive integer, and $T$ a topological $n$-ball. A \emph{topoplane} in $T$ is a topological $(n-1)$-ball $H \subseteq T$ that divides $T$ into two connected topological subspaces.
\end{definition}

\begin{definition}
Let $\mathscr{A}$ be a finite set of topoplanes in a topological $n$-ball $T$. A \emph{flat} of $\mathscr{A}$ is a nonempty intersection of topoplanes in $\mathscr{A}$. Denote by $L(\mathscr{A})$ the set composed by the flats of $\mathscr{A}$. 
\end{definition}

\begin{example}
The flat set generated by both topoplanes in the yellow open disk of Figure~\ref{Pic1} is composed of the yellow disk, both topoplanes, and the four intersection points.
	
\begin{figure}[h]
\centering
\includegraphics[scale=0.17]{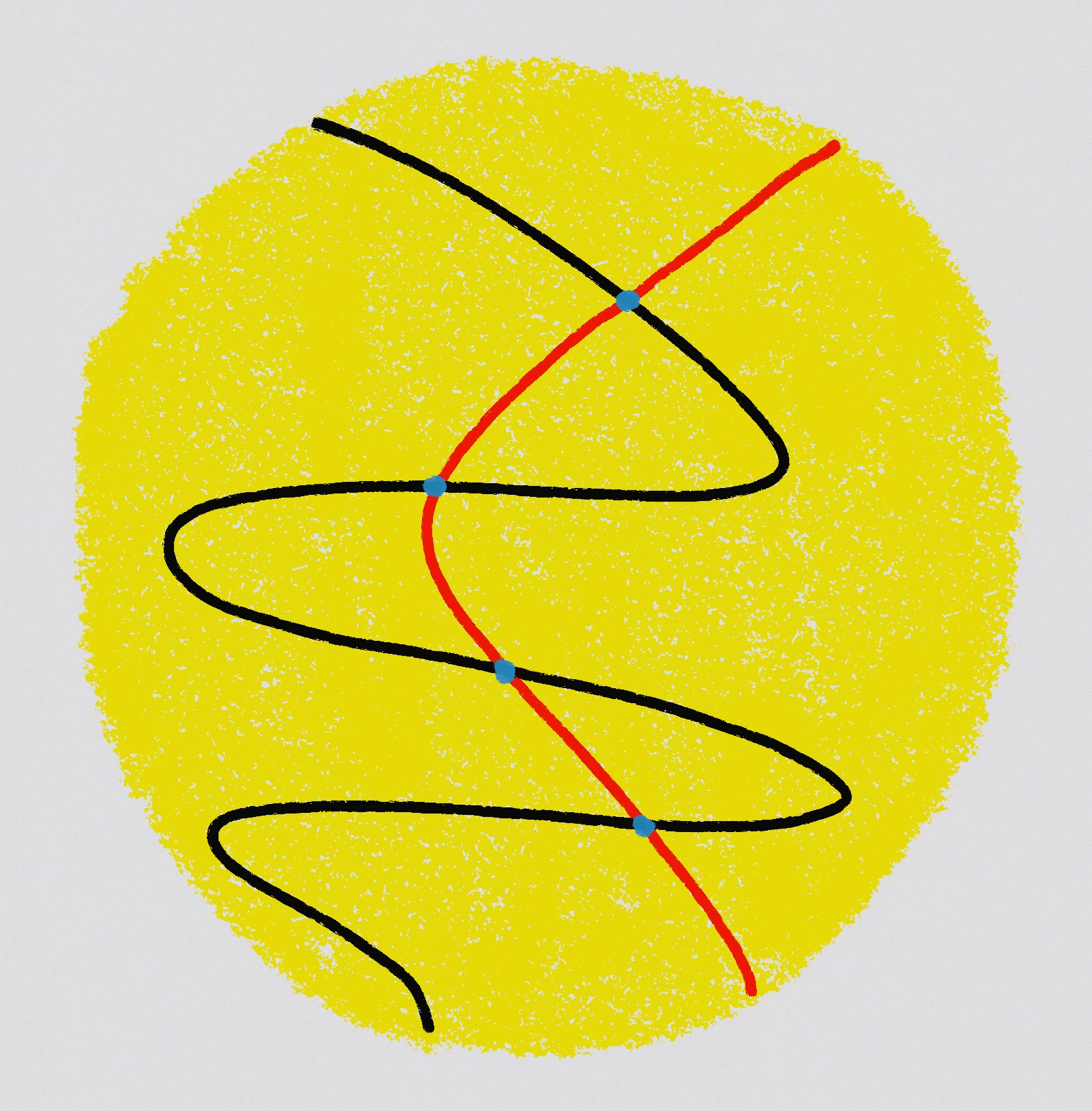}
\caption{Two Topoplanes in an Open Disk}
\label{Pic1}
\end{figure}
\end{example}

\begin{definition} \label{DefArr}
Let $\mathscr{A}$ be a finite set of topoplanes in a topological ball $T$. It is a \emph{topoplane arrangement} if
\begin{itemize}
\item[(i)] every flat in $L(\mathscr{A})$ is a topological ball,
\item[(ii)] for every topoplane $H \in \mathscr{A}$ and each flat $X \in L(\mathscr{A})$, either $X \subseteq H$ or $H \cap X = \varnothing$ or $H \cap X$ is a topoplane in $X$.
\end{itemize}
\end{definition}

\begin{example}
The flat set of the topoplane arrangement in Figure~\ref{Pic2} is composed of $\mathbb{R}^3$, both topoplanes, and the intersection point. 
	
\begin{figure}[h]
\centering
\includegraphics[scale=0.19]{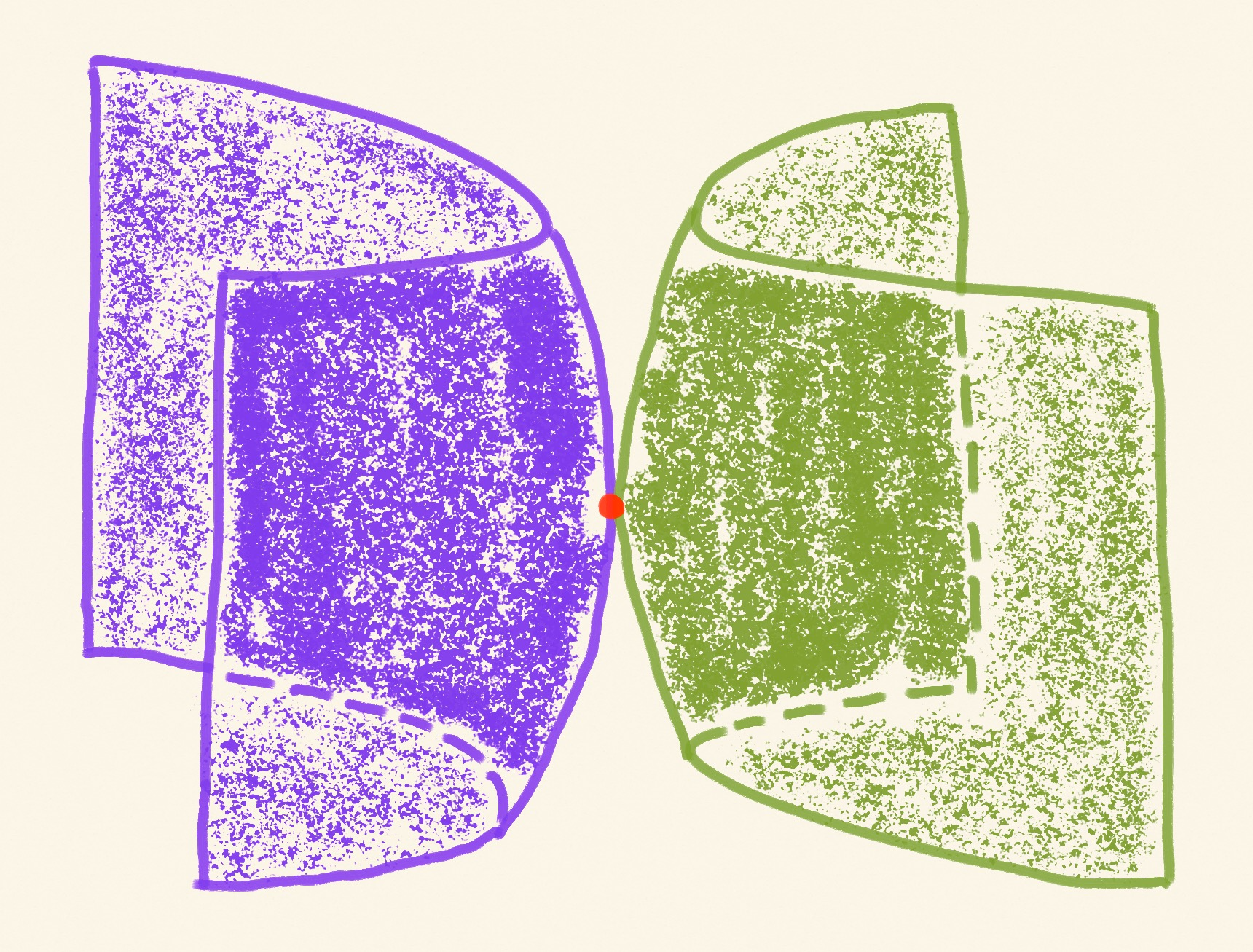}
\caption{A Topoplane Arrangement in $\mathbb{R}^3$}
\label{Pic2}
\end{figure}
\end{example}

\begin{proposition}\citep[Prop.~1]{forge2009division} \label{PrRes}
Let $\mathscr{A}$ be a topoplane arrangement in a topological ball $T$, and consider a flat $X \in L(\mathscr{A})$. The induced set of topological subspaces in $X$ defined by
$$\mathscr{A}^X := \{X \cap H\ |\ H \in \mathcal{A},\, X \nsubseteq H,\, X \cap H \neq \varnothing\}$$ is a topoplane arrangement in $X$.
\end{proposition}

\begin{definition}
Let $\mathscr{A}$ be a topoplane arrangement in a topological ball $T$, and consider a flat $X \in L(\mathscr{A})$. The topoplane arrangement $\mathscr{A}^X := \{X \cap H\ |\ H \in \mathscr{A},\, X \nsubseteq H,\, X \cap H \neq \varnothing\}$ in $X$ is called the \emph{restriction} of $\mathscr{A}$ on $X$.
\end{definition}

\begin{definition}
Let $\mathscr{A}$ be a finite set of topoplanes in a topological ball $T$. A pair of distinct topoplanes $(H,\,K) \in \mathscr{A} \times \mathscr{A}$ forms a \emph{transsection} if $H \setminus K$ is composed by two components which lie on opposite sides of $K$.
\end{definition}

\begin{definition}
Let $\mathscr{A}$ be a topoplane arrangement in a topological ball $T$. It is said to be \emph{transsective} if, for each pair of distinct topoplanes $(H,\,K) \in \mathscr{A} \times \mathscr{A}$, either $H \cap K = \varnothing$ or $(H,\,K)$ forms a transsection.
\end{definition}

\begin{proposition} \label{PrTra}
Let $\mathscr{A}$ be a topoplane arrangement in a topological ball $T$, and $X$ a flat in $L(\mathscr{A})$. If $\mathcal{A}$ is transsective, then $\mathscr{A}^X$ is a transsective topoplane arrangement in $X$. 
\end{proposition}

\begin{proof}
Consider two distinct topoplanes in $X$, namely having the forms $X \cap H$ and $X \cap K$ with $H, K \in \mathscr{A}$. Suppose that $(X \cap H) \cap (X \cap K) \neq \varnothing$. Since $X \cap H \nsubseteq X \cap K$ and $\mathscr{A}^X$ is a topoplane arrangement as seen in Proposition~\ref{PrRes}, then $(X \cap H) \cap (X \cap K) = X \cap H \cap K$ is a topoplane in $X \cap H$. Hence $X \cap H \cap K$ divides $X \cap H$ into two connected topological subspaces $(X \cap H \cap K)^1$ and $(X \cap H \cap K)^{-1}$. Besides, the topoplane $X \cap K$ divides $X$ into two connected topological subspaces $(X \cap K)^1$ and $(X \cap K)^{-1}$, and we have
\begin{itemize}
\item either $(X \cap H \cap K)^1 \subseteq (X \cap K)^1$ and $(X \cap H \cap K)^{-1} \subseteq (X \cap K)^{-1}$,
\item or $(X \cap H \cap K)^{-1} \subseteq (X \cap K)^1$ and $(X \cap H \cap K)^1 \subseteq (X \cap K)^{-1}$.
\end{itemize}
In both cases, $(X \cap H) \setminus (X \cap K)$ is composed by two components in $X$ which lie on opposite sides of $X \cap K$. The topoplane arrangement $\mathscr{A}^X$ is consequently transsective.
\end{proof}

\begin{definition}
Let $\mathscr{A}$ be a transsective topoplane arrangement in a topological ball $T$. Denote by $H^{-1}$ and $H^1$ both connected components obtained after division of $T$ by a topoplane $H \in \mathscr{A}$. Moreover, set $H^0 = H$. The \emph{sign map} of $H$ is the function
$$\sigma_H: T \rightarrow \{-1,\,0,\,1\}, \quad v \mapsto \begin{cases}
-1 & \text{if}\ v \in H^{-1}, \\
0 & \text{if}\ v \in H^0, \\
1 & \text{if}\ v \in H^1. 
\end{cases}$$
The sign map of $\mathscr{A}$ is the function $\sigma_{\mathscr{A}}: T \rightarrow \{-1,\,0,\,1\}^{\mathscr{A}}, \ v \mapsto \big(\sigma_H(v)\big)_{H \in \mathscr{A}}$. And the \emph{sign set} of $\mathscr{A}$ is the set $$\sigma_{\mathscr{A}}(T) := \big\{\sigma_{\mathscr{A}}(v)\ \big|\ v \in T\big\}.$$
\end{definition}

\begin{definition}
Let $\mathscr{A}$ be a transsective topoplane arrangement in a topological ball $T$. A \emph{face} of $\mathscr{A}$ is a subset $F$ of $T$ such that
$$\exists x \in \sigma_{\mathscr{A}}(T),\ F = \big\{v \in T\ \big|\ \sigma_{\mathscr{A}}(v) = x\big\}.$$
A \emph{chamber} of $\mathscr{A}$ is a face $F$ such that $\sigma_{\mathscr{A}}(F) \in \{-1,\,1\}^{\mathscr{A}}$. Denote by $F(\mathscr{A})$ and $C(\mathscr{A})$ the sets composed by the faces and the chambers of $\mathscr{A}$, respectively. 
\end{definition}

\begin{proposition} \label{PrDef}
Let $\mathscr{A}$ be a transsective topoplane arrangement in a topological ball $T$. Then, every face of $\mathscr{A}$ is a deformed ball.
\end{proposition}

\begin{proof}
Assume $T$ is a topological $n$-ball, and begin by considering a chamber $C \in C(\mathscr{A})$:
\begin{itemize}
\item Let $x \in C$, and $d = \min \big\{\mathrm{dist}(x, H)\ \|\ H \in \mathscr{A}\big\}$, where $\mathrm{dist}$ is a distance function on $T$. Then, the $n$-ball of radius $d/2$ and center $x$ is included in $C$. The chamber $C$ is consequently an $n$-manifold.
\item Let $x,y \in C$. The fact that $\mathscr{A}$ is transsective and $\sigma_{\mathscr{A}}(x) = \sigma_{\mathscr{A}}(y)$ imply the path connectivity of $x$ and $y$.
\item The chamber $C$ can naturally not contain holes, meaning that $\pi_k(C, x_0)$ is trivial for each positive integer $k$ and distinguished point $x_0$ of $C$.
\end{itemize} 
The chamber $C$ is then a deformed ball. Consider now a face $F \in F(\mathscr{A}) \setminus C(\mathscr{A})$, and the flat $\displaystyle X = \bigcap_{\substack{H \in \mathscr{A} \\ \sigma_H(F) = 0}} H$. We know from Proposition~\ref{PrTra} that $\mathscr{A}^X$ is a transsective topoplane arrangement in $X$. As $F$ is a chamber of $\mathscr{A}^X$, it is therefore a deformed ball.
\end{proof}

\section{The $f$-Polynomial of a Topoplane Arrangement} \label{fPo}

\noindent We finally get the $f$-polynomial of a transsective topoplane arrangement $\mathscr{A}$ in a topological ball $T$ in Theorem~\ref{fPol} of this section. Besides, investigating the constant of that polynomial gives that $\displaystyle \sum_{X \in L(\mathscr{A})} (-1)^{\mathrm{rank}\,X} \mu(T,X)$ is the number of chambers of $\mathscr{A}$.

\begin{definition}
Let $\mathscr{A}$ be a transsective topoplane arrangement in a topological ball. Define the \emph{dimension} $\dim X$ of a flat $X$ of $\mathscr{A}$ which is topological $n$-ball, as well as the dimension $\dim F$ of a face $F$ of $\mathscr{A}$ which is a deformed $n$-ball, to be $n$. Call such flat and face of $\mathscr{A}$ $n$-flat and $n$-face, respectively.
\end{definition}

\begin{definition}
Consider a transsective topoplane arrangement $\mathscr{A}$ in a topological $n$-ball. Let $f_i(\mathscr{A})$ be the number of $i$-faces of $\mathscr{A}$, and $x$ a variable. The \emph{$f$-polynomial} of $\mathscr{A}$ is
$$f_{\mathscr{A}}(x) := \sum_{i=0}^n f_i(\mathscr{A}) \, x^{n-i}.$$
\end{definition}

\begin{definition}
Let $\mathscr{A}$ be a transsective topoplane arrangement in a topological $n$-ball. Define the \emph{rank} of a flat $X \in L(\mathscr{A})$ to be $\mathrm{rank}\,X := n - \dim X$, and that of the topoplane arrangement $\mathscr{A}$ to be $\mathrm{rank}\,\mathscr{A} := \max \big\{\mathrm{rank}\,X \in \mathbb{N}_0\ \big|\ X \in L(\mathscr{A})\big\}$. 
\end{definition}

\noindent Recall that the Möbius function $\mu: L(\mathscr{A}) \times L(\mathscr{A}) \rightarrow \mathbb{Z}$ of a meet semilattice $L(\mathscr{A})$ is recursively defined, for $X,Y \in L(\mathscr{A})$, by $$\displaystyle \mu(X,Y) := \begin{cases}
1 & \text{if}\ X = Y, \\
\displaystyle -\sum_{\substack{Z \in L(\mathscr{A}) \\ X \leq Z < Y}} \mu(X,Z) = -\sum_{\substack{Z \in L(\mathscr{A}) \\ X < Z \leq Y}} \mu(Z,Y) & \text{if}\ X < Y, \\
0 & \text{otherwise} 
\end{cases}.$$

\begin{definition}
Let $\mathscr{A}$ be a transsective topoplane arrangement in a topological ball, and $x,y$ variables. The \emph{Möbius polynomial} of $\mathscr{A}$ is $$M_{\mathscr{A}}(x,y) := \sum_{X,Y \in L(\mathscr{A})} \mu(X,Y) \, x^{\mathrm{rank}\,X} \, y^{\mathrm{rank}\,\mathscr{A} - \mathrm{rank}\,Y}.$$
\end{definition}

\noindent We can now state the main result of this article.

\begin{theorem} \label{fPol}
Let $\mathscr{A}$ be a transsective topoplane arrangement in a topological ball. The $f$-polynomial of $\mathscr{A}$ is $$f_{\mathscr{A}}(x) = (-1)^{\mathrm{rank}\,\mathscr{A}} M_{\mathscr{A}}(-x,-1).$$
\end{theorem}

\begin{proof}
We know from Proposition~\ref{PrTra} and Proposition~\ref{PrDef} that the pair $\big(X,\, F(\mathscr{A}^X)\big)$ forms a system of deformed balls. Thus
$$\displaystyle \chi(X) = \sum_{i = 0}^{\dim X} (-1)^i f_i(\mathscr{A}^X) = (-1)^{\dim X}.$$
Moreover, as every $i$-face $F \in F(\mathscr{A}^X)$ is a chamber of a unique $i$-flat $\displaystyle \bigcap_{\substack{H \in \mathscr{A}^X \\ \sigma_H(F) = 0}} H \in L(\mathscr{A}^X)$, then $\displaystyle f_i(\mathscr{A}^X) = \sum_{\substack{Y \in L(\mathscr{A}^X) \\ \dim Y = i}} \#C\Big({(\mathscr{A}^X)}^Y\Big)$, and $$\sum_{Y \in L(\mathscr{A}^X)} (-1)^{\dim Y} \#C\Big({(\mathscr{A}^X)}^Y\Big) = (-1)^{\dim X}.$$
We have $L(\mathscr{A}^X) = \big\{Y \in L(\mathscr{A})\ \big|\ Y \geq X\big\}$, and, for every $Y \in L(\mathscr{A}^X)$, also $C\Big({(\mathscr{A}^X)}^Y\Big) = C(\mathscr{A}^Y)$. Hence,
$$\sum_{\substack{Y \in L(\mathscr{A}) \\ Y \geq X}} (-1)^{\dim Y} \#C(\mathscr{A}^Y) = (-1)^{\dim X}.$$
Using the Möbius inversion formula, we obtain $$\sum_{\substack{Y \in L(\mathscr{A}) \\ Y \geq X}} (-1)^{\dim Y} \mu(X,Y) = (-1)^{\dim X} \#C(\mathscr{A}^X).$$
Besides, $\displaystyle (-1)^{\mathrm{rank}\,\mathscr{A}} M_{\mathscr{A}}(-x,-1) = \sum_{X,Y \in L(\mathscr{A})} (-1)^{\dim Y - \dim X} \mu(X,Y)\, x^{\mathrm{rank}\,X}$. Therefore, for every $0 \leq i \leq n$, the coefficient $\lambda_{n-i}$ of $x^{n-i}$ in the polynomial $(-1)^{\mathrm{rank}\,\mathscr{A}} M_{\mathscr{A}}(-x,-1)$ is $$\lambda_{n-i} = \sum_{\substack{X \in L(\mathscr{A}) \\ \dim X = i}} \sum_{Y \in L(\mathscr{A}^X)} (-1)^{\dim Y - \dim X} \mu(X,Y) = \sum_{\substack{X \in L(\mathscr{A}) \\ \dim X = i}} \#C(\mathscr{A}^X) = f_i(\mathscr{A}).$$
\end{proof}

\begin{example}
Consider the topoplane arrangement $\mathscr{A}_{\text{ex}}$ formed by nine topoplanes in $\mathbb{R}^2$ represented in Figure~\ref{Ex2}. As its Möbius polynomial is $M_{\mathscr{A}_{\text{ex}}}(x,y) = 5x^2 + y^2 + 9xy - 11x - 9y + 6$, its $f$-polynomial is then $f_{\mathscr{A}_{ex}}(x) = 5x^2 + 20x + 16$.

\bigskip

\begin{figure}[h]
\centering
\includegraphics[scale=0.19]{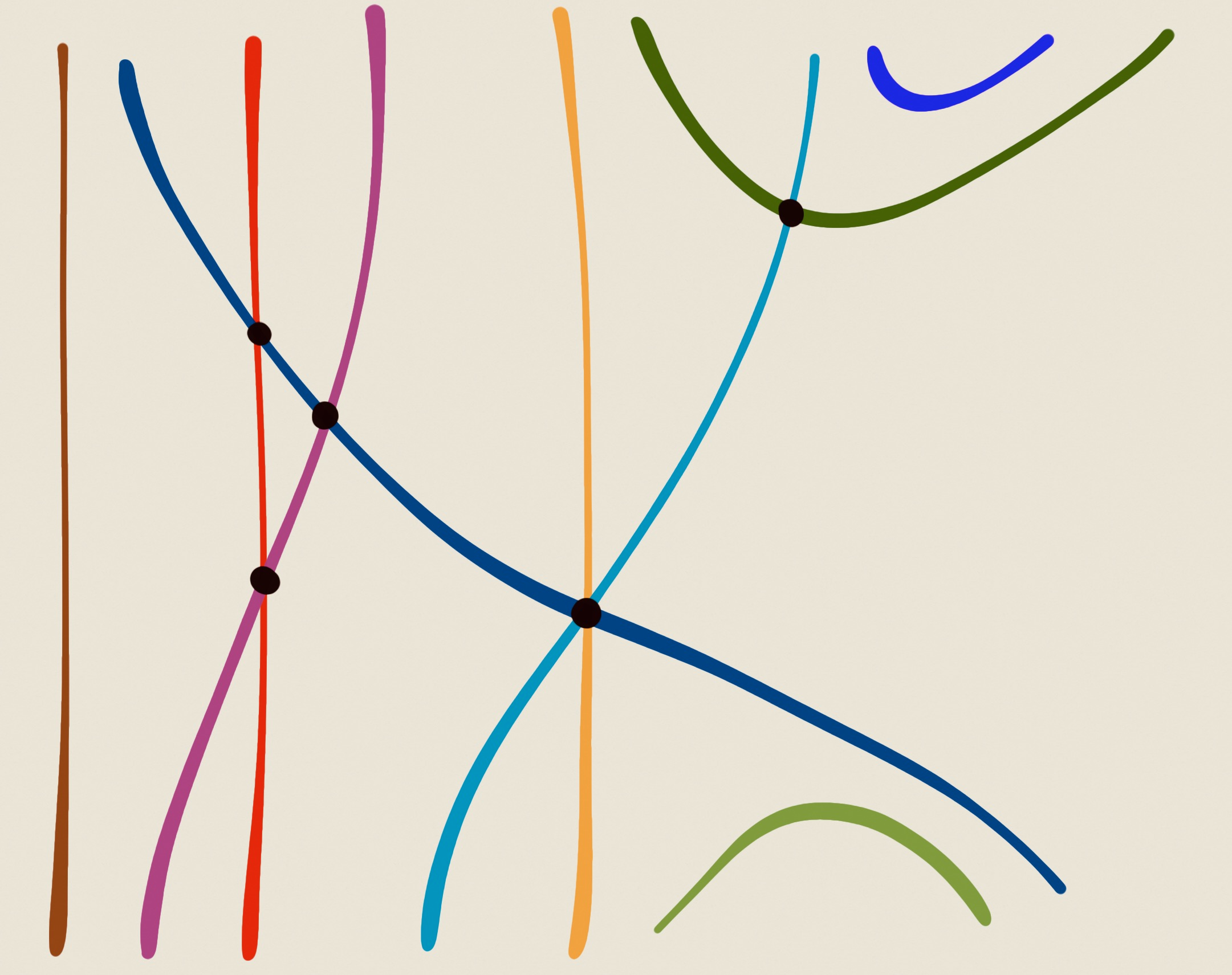}
\caption{The Topoplane Arrangement $\mathscr{A}_{\text{ex}}$}
\label{Ex2}
\end{figure}
\end{example}

\bibliographystyle{agsm}
\bibliography{References}

\end{document}